\documentclass[12pt,a4paper,reqno]{amsart}
\usepackage{amssymb}
\usepackage{latexsym}
\usepackage{exscale}
\usepackage{epsfig}
\usepackage{verbatim}
\usepackage{float,graphicx}
\usepackage{subfigure}
\usepackage{epic,eepic}
\usepackage[latin1]{inputenc}  

\newcommand{\abs}[1]{\left| #1 \right|}
\newcommand{\norm}[1]{\| #1 \|}

\newcommand{\ip}[2]{\langle #1 , #2 \rangle}

\numberwithin{equation}{section}
\newtheorem{theorem}{Theorem}[section]

\begin{document}

\allowdisplaybreaks
\title[A simple shearlet-based reconstruction for CT]  
{\large A simple shearlet-based reconstruction for computer tomography} 

\vskip 1cm
\author{Santiago Córdova}
\address{Santiago Córdova
\\
Departamento de Matemáticas
\\
ITAM (Instituto Tecnológico Autónomo de México)
\\
Río Hondo No. 1, Col. Tizapán CP 01080 Ciudad de México, México}

\email{santiago.cordova@itam.mx}

\vskip 1cm
\author{Daniel Vera}
\address{Daniel Vera
\\
Departamento de Matemáticas
\\
ITAM (Instituto Tecnológico Autónomo de México)
\\
Río Hondo No. 1, Col. Tizapán CP 01080 Ciudad de México, México}

\email{daniel.vera@itam.mx}


\date{\today}

\maketitle

\begin{abstract} We find a new and simple inversion formula of the Radon transform RT
with the only use of the shearlet system and of well-known properties of RT. No intertwining relation of differential operators in Euclidean space and Radon domain is used. As a consequence, an additive noise is not incremented. Since the continuum theory of shearlets has a straight translation to the discrete theory, we find a fast, stable and computable algorithm that recovers a digital image from noisy samples of the Radon transform preserving edges. In the process, we find a more natural and easier-to-construct ``density-compensation weight" functions for the ShearLab toolbox.
\end{abstract}

\section{Introduction}\label{S:Intro}
In computer tomography the problem is to reconstruct an image from noisy projections of X-rays. This inverse problem is a non-invasive technique used not only in medicine but also in industrial scanning. It is related to other inverse problems as in seismology and MRI, for example. CT scanning in 2D is intimately related to the Radon transform in 2D, which is the problem of recovering a density function from line integrals. This problem was solved in 1917 by Johann Radon in \cite{Rad17} when the whole data is known, i.e. for uncountable angles as well as distances, as follows. Let $\Theta=(\cos\theta, \sin\theta)$ be a vector in the unitary sphere $S^1$. For $x\in\mathbb{R}^2$ and $s\in \mathbb{R}$ define formally the Radon transform as line integrals with direction $\Theta$ and distance $s$
\begin{equation}\label{e:def_RadonTrnsf}
  \mathcal{R}f(\theta,s):=\int_{\ip{x}{\Theta}=s} f(x)dx
        =\int_{\mathbb{R}^2} f(x)\delta(s-\ip{x}{\Theta}) dx,
\end{equation}
where $\delta$ is Dirac's distribution. For $g(\theta,s)=\mathcal{R}f(\theta,s)$, Radon's inversion formula formally reads as
\begin{equation}\label{e:Radon-Inversion}
  f(x)=\frac{1}{4\pi^2}\int_{S^1}\int_{\mathbb{R}} \frac{\frac{d}{ds}g(\theta,s)}{\ip{x}{\Theta}-s}dsd\Theta.
\end{equation}\vskip0.2cm
How to implement the inversion formula (\ref{e:Radon-Inversion}) with accurate, efficient, fast and computable algorithms is far from obvious. Starting with the problem of measuring
the projections, several acquisition methods have been proposed, v.gr. parallel-beam or fan-beam in the 2D problem. Whereas for the solution, one can find in the broad literature that there exist also several solving methods like the filtered back-projection (FBP), iterative or algebraic methods. It is well-known that the FBP is equivalent to the Radon inversion formula (see\cite{NaWu00}).

More recently, multi-resolution analysis MRA methods have succeeded in showing their very high performance when compared to more traditional methods. We mention only two inspiring works, that of C. Berenstein and D. Walnut in \cite{BeWa94} and that of S. Bonnet, F. Peyrin, F. Turjman and R. Prost in \cite{BPTP02}, and references therein. Nevertheless, these methods try to find good discretizations of an integral relating wavelet and Radon transforms in the time domain. They are, therefore, approximation algorithms. These methods are based on the FBP method and, as a consequence, they need to smooth the unbounded filtering operator in the time domain.

Since then there have been advances in two directions. First, the development of the shearlet system which is a
multi-scale and multi-directional system defined by affine operators over a single function that act on the integer grid only (see \cite{GKL06}, for example) and have a tiling of the frequency plane with resemblance to the wavelet induced tiling. Second, the pseudo-polar fast Fourier transform (or PP-FFT for short) as developed by A. Averbuch, R. Coifman, D. L. Donoho and M. Israeli in \cite{ACDIS} which will allow us to discretize the main result. This fundamental article also define a discrete Radon transform with a fast and exact inverse, however it lacks an analysis on noisy data and does not have the advantages of multi-scale methods.

Shearlets are just an example of a variety of multi-scale and multi-directional systems like curvelets \cite{CaDo02} and contourlets \cite{DoVe}. Unlike curvelets, which are defined in polar coordinates, shearlets are defined by affine transformations which make them more appealing to computable procedures. And, unlike contourlets which are defined by digital filters, there exist a continuous theory of shearlets which naturally translates to a discrete theory.

E. J. Candès and D. L. Donoho applied the curvelets to recover a density function from noisy Radon data in the breakthrough paper \cite{CaDo02}. This was achieved by applying operators of fractional differentiation to the curvelet system and a well-known intertwining relation between fractional powers of the Laplacian in $\mathbb{R}^2$ and fractional differentiation along the $s$-coordinate in the Radon domain. The drawback is that the reproducing formula (by the quasi-biorthogonal systems) of $f$ has a factor of $2^j, j\geq 0$, in the coefficients at scale $2^{-j}$, making the inverse problem highly ill-posed because noise is a high frequency phenomena. Another drawback is that when white noise model is used on the Radon domain, that is, independent identically distributed i.i.d. random variables $N(0,\sigma^2)$, the noise is not i.i.d. in the representation of $f$ by the Radon-transformed curvelet system since this is not an orthonormal transformation. Estimation of $f$ is done by a thresholding dependent of the scale.

In \cite{CEGL}, F. Colonna, G. Easly, K. Guo and D. Labate, obtained the inversion of the Radon transformed applying the shearlet system in the same way as in \cite{CaDo02} with the same problems and solutions.

The main result in this paper does not use the intertwining relation to evaluate the inverse Radon transform. One of the remarkable features is that, since we use the Fourier slice Theorem before evaluation of shearlet coefficients, an i.i.d. $N(0,\sigma^2)$ vector in the Radon data is mapped to an i.i.d. $N(0,\sigma^2)$ vector in the PP Fourier data since the Fourier transform is orthonormal. The noise is then homoscedastic and, therefore, thresholding does not depend on scales. Noise is not magnified in high-frequencies by a factor $2^j$ as in \cite{CaDo02} and \cite{CEGL}. Unlike \cite{CaDo02} and \cite{CEGL}, where the coarse-scale index of the curvelet/shearlet system should be set according to the noise level (usually it runs from $j\geq 0$), our starting index need not be dependent of the noise level.

It will be obvious from our results that a similar technique can be used to prove analog results in the framework of curvelets.

In a more recent paper \cite{BdMdVO}, F. Bartolucci, F. De Mari, E. De Vito and F. Odone, find a similar result to ours. We point out some differences between \cite{BdMdVO} and our results: i) \cite{BdMdVO} is developed for continuous parameters of scale, shearing, and location, ii) evaluation of the shearlet transform is done in time domain, iii) \cite{BdMdVO} admits that the results only ``open new perspectives in the inversion of the Radon transform", iv) a discussion on noisy data is missing, and v) \cite{BdMdVO} uses, we believe, very sophisticated machinery.

Acknowledgements. Theorem \ref{t:RadonInv-Shrlts} was obtained when the second author was doing a postdoctorate in the University of Houston. He thanks Demetrio Labate and Bart Goossens for very enlightening talks on CT.

\section{Background}
We will first present a series of well-known properties of the Radon transform in 2D. They are stated for functions $f\in\mathcal{S}$, however the results can be extended to other function spaces, see \cite{BeWa94} and references therein. Afterwards, we introduce briefly the shearlet system. We close this section defining the linogram acquisition method and its relation with the pseudo-polar grid.
\subsection{Radon transform and properties}\label{sS:RT}
We follow closely section 2 in Natterer and Wübbeling's book \cite{NaWu00}. Let $\Theta=(\cos\theta,\sin\theta)$ and $x\in\mathbb{R}^2$. The Radon transform has been defined in
(\ref{e:def_RadonTrnsf}). Given $g(\theta,s)=\mathcal{R}f(\theta, s)$, define the convolution and Fourier transform in the second variable as
\begin{equation}\label{e:Conv_2ndVar}
    (g\ast h)(\theta,s):=\int_\mathbb{R} g(\theta, s-t)h(\theta, t)dt,
\end{equation}
and
\begin{equation}\label{e:FT_2ndVar}
    \hat{g}(\theta,\sigma):=\int_\mathbb{R} g(\theta,s)e^{-2\pi i\sigma s} ds=\mathcal{F}_{s\rightarrow \sigma}(g(\theta, \cdot))(\sigma),
\end{equation}
respectively. 
The following is known as the ``projection theorem", ``central slice theorem" or as the ``Fourier slice"
\begin{theorem}\label{T:Fourier-slice}
    Let $f\in \mathcal{S}(\mathbb{R}^2)$. Then, for fix $\Theta\in S^1$ and $\sigma\in \mathbb{R}$,
    \begin{equation}\label{e:Fourier-slice}
        (\mathcal{R}f)^\wedge(\theta,\sigma)=\hat{f}(\sigma\Theta).
    \end{equation}
\end{theorem}This means that, for a fixed $\theta$, the Fourier transform of $\mathcal{R}f$ is equal to the usual 2D Fourier transform of $f$
restricted to the slice $\sigma\Theta$. Another useful result is the next
\begin{theorem}\label{T:Conv-equiv}
    Let $f,g\in \mathcal{S}(\mathbb{R}^2)$. Then,
        \begin{equation}\label{e:Conv-equiv}
            \mathcal{R}f\ast \mathcal{R}g=\mathcal{R}(f\ast g).
        \end{equation}
\end{theorem}
Let us introduce the backprojection operator for $g\in \mathcal{S}(C^2)$, where $C^2=S^1\times \mathbb{R}$ is sometimes known as the unitary cilinder.
\begin{equation}\label{e:BP-op}
    \mathcal{R}^\ast g(x):=\int_0^{2\pi} g(\theta,\ip{x}{\Theta}) d\theta.
\end{equation}
For $g=\mathcal{R}f$, $\mathcal{R}^\ast g$ is the average of the line integrals of $f$ through $x$, $\mathcal{R}^\ast $ is indeed the adjoint of $\mathcal{R} $.
We also rely our results on the next
\begin{theorem}\label{T:BP-equiv}
    Let $f\in\mathcal{S}(\mathbb{R}^2)$ and $g\in \mathcal{S}(C^2)$. Then,
    \begin{equation}\label{e:BP-equiv}
        (\mathcal{R}^\ast g)\ast f = \mathcal{R}^\ast(g\ast \mathcal{R}f).
    \end{equation}

\end{theorem}
A whole family of exact inversion formulas uses the Riesz potentials in $\mathbb{R}^d$ or $C^d$ defined as
$(I^\alpha)^\wedge(\xi):=\abs{\xi}^{-\alpha}\hat{f}(\xi), \xi\in\mathbb{R}^d$ with $\alpha<d$, and
$(I^\alpha g)^\wedge(\theta,\sigma):=\abs{\sigma}^{-\alpha}\hat{g}(\theta,\sigma), \alpha<1$, respectively.
\begin{theorem}\label{T:InvFrmls-RieszPotencials}
    Let $f\in\mathcal{S}(\mathbb{R}^d)$ and $g=\mathcal{R}f$. Then, for $\alpha<d$
    \begin{equation}\label{e:InvFrmls-RieszPotencials}
        f=\frac{1}{2}I^{-\alpha}\mathcal{R}^\ast I^{\alpha-d+1} g.
    \end{equation}
\end{theorem}
As a special case we have for $\alpha=0$ and dimension $d=2$,
\begin{equation}\label{e:InvFrmls-RieszPotencials-special-case}
f=\frac{1}{2} \mathcal{R}^\ast I^{-1} g = \frac{1}{2} \mathcal{R}^\ast I^{-1} \mathcal{R} f,
\end{equation}
where
$(I^{-1}g)^\wedge(\theta,\sigma)=\abs{\sigma}\hat{g}(\theta,\sigma)=\abs{\sigma}\hat{f}(\sigma\Theta)$ because of the Fourier slice Theorem.
The factor $\abs{\sigma}$ is a Fourier multiplier (or a filter) which is unbounded.

\vskip1cm

\subsection{Shearlets}\label{S:Shrlts}
Define the horizontal cone
\begin{equation*}\label{e:Cone_Domain}
\mathcal{D}^h:=\{(\xi_1,\xi_2)\in\hat{\mathbb{R}}^2:
    \abs{\xi_1}\geq\frac{1}{8}, \abs{\frac{\xi_2}{\xi_1}}\leq 1\}.
\end{equation*}
Let now $\hat{\psi}_1\in C^\infty(\mathbb{R})$, be a wavelet with
$\text{supp }\hat{\psi}_1\subset
[-\frac{1}{2},-\frac{1}{16}]\cup[\frac{1}{16},\frac{1}{2}]$ such that
\begin{equation*}\label{e:Discrt_Shrlt_Cond_Cone_1}
\sum_{j\geq 0}\abs{\hat{\psi}_1(2^{-2j}\omega)}^2 =1, \;\;\;
\text{for }\abs{\omega}\geq \frac{1}{8},
\end{equation*}
and $\hat{\psi}_2\in C^\infty(\mathbb{R})$ be a bump function with
$\text{supp }\hat{\psi}_2\subset [-1,1]$ such that
\begin{equation*}\label{e:Discrt_Shrlt_Cond_Cone_2}
\abs{\hat{\psi}_2(\omega-1)}^2+\abs{\hat{\psi}_2(\omega)}^2+\abs{\hat{\psi}_2(\omega+1)}^2=1,
\;\;\; \text{for } \abs{\omega}\leq 1.
\end{equation*}
It follows that, for
$j\geq 0$,
\begin{equation*}\label{e:Discrt_Shrlt_Cond_Cone_3}
\sum_{\ell=-2^j}^{2^j} \abs{\hat{\psi}_2(2^j\omega-\ell)}^2=1,
\;\;\; \text{for }\abs{\omega}\leq 1.
\end{equation*}
Now let
$$A_h= \left(%
\begin{array}{cc}
  4 & 0 \\
  0 & 2 \\
\end{array}%
\right), \;\;\; B_h=\left(%
\begin{array}{cc}
  1 & 1 \\
  0 & 1 \\
\end{array}%
\right)$$ and
$\hat{\psi}^h(\xi)=\hat{\psi}_1(\xi_1)\hat{\psi}_2(\frac{\xi_2}{\xi_1})$.
It follows that
\begin{eqnarray*}\label{e:PrsvlFrm_Prop_Shrlts}
\nonumber
  \sum_{j\geq 0}\sum_{\ell=-2^j}^{2^j} \abs{\hat{\psi}^h(\xi A^{-j}_hB^{-\ell}_h)}^2
    &=& \sum_{j\geq 0}\sum_{\ell=-2^j}^{2^j} \abs{\hat{\psi}_1(2^{-2j}\xi_1)}^2
        \abs{\hat{\psi}_2(2^j\frac{\xi_2}{\xi_1}-\ell)}^2 \\
    &=& 1
\end{eqnarray*}
for $\xi=(\xi_1,\xi_2)\in\mathcal{D}^h$. This is the
Parseval frame condition for the horizontal cone. Since
$\text{supp }\hat{\psi}^h\subset [-\frac{1}{2},\frac{1}{2}]^2$,
(\ref{e:PrsvlFrm_Prop_Shrlts}) implies that the shearlet system
\begin{equation*}\label{e:Shrlt_Sys_Cone}
\{\psi_{j,\ell,k}^h(x)= 2^{3j/2}\psi^h(B^\ell_h A^j_h x-k): j\geq
0, -2^j\leq\ell\leq2^j, k\in\mathbb{Z}^2\},
\end{equation*}
is a Parseval frame for $L^2((\mathcal{D}^h)^\vee)=\{f\in
L^2(\mathbb{R}^2): \text{supp }\hat{f}\subset\mathcal{D}^h\}$, this means that
$$\sum_{j\geq 0}\sum_{\ell=-2^j}^{2^j}\sum_{k\in\mathbb{Z}^2}
    \abs{\ip{f}{\psi_{j,\ell,k}^h}}^2 = \norm{f}^2_{L^2(\mathbb{R}^2)},$$
for all $f\in L^2(\mathbb{R}^2)$ such that $\text{supp
}\hat{f}\subset \mathcal{D}^h$. One can also construct a Parseval
frame for the vertical cone $\mathcal{D}^v$ and add a low-frequency function so that the complete shearlet system is a Parseval frame for the whole space $L^2(\mathbb{R}^2)$
as follows. Let $\hat{\varphi}\in C^\infty_c(\mathbb{R}^2)$, with $\text{supp
}\hat{\varphi}\subset [-\frac{1}{4},\frac{1}{4}]^2$ and
$\abs{\hat{\varphi}}=1$ for $\xi\in
[-\frac{1}{8},\frac{1}{8}]^2=\mathcal{R}$, be such that
\begin{eqnarray*}\label{e:Shrlt_Repr_Sys}
\nonumber
1
    &=& \abs{\hat{\varphi}(\xi)}^2\chi_\mathcal{R}(\xi)
    + \sum_{j\geq 0}\sum_{\ell=-2^j}^{2^j}\abs{\hat{\psi}^h(\xi A^{-j}_hB^{-\ell}_h)}^2\chi_{\mathcal{D}^h}(\xi) \\
    &&\;\;\;\; +\sum_{j\geq 0}\sum_{\ell=-2^j}^{2^j}\abs{\hat{\psi}^v(\xi
        A^{-j}_vB^{-\ell}_v)}^2\chi_{\mathcal{D}^v}(\xi),
        \;\;\;\text{for all } \xi\in\hat{\mathbb{R}}^2.
\end{eqnarray*}
For a more comprehensive introduction see \cite{GKL06}. Let
\begin{equation}\label{e:Shrlt_Sys_Cone}
\{\psi_{j,\ell,k}^{(\mathfrak{d})}(x)= \abs{\text{det
}A_{(\mathfrak{d})}}^{j/2}\psi^{(\mathfrak{d})}(B^{\ell}_{(\mathfrak{d})}
A^j_{(\mathfrak{d})} x-k): \mathfrak{d}=\{h,v\}, j\geq 0,
\abs{\ell}\preceq2^j, k\in\mathbb{Z}^d\},
\end{equation}
be the shearlet system associated to the horizontal $h$ or vertical
$v$ cone of frequencies. Let $P=(\mathfrak{d},j,\ell,k)$ be the
indexing associated to the parallelograms
\begin{equation}\label{e:def_Qjlk}
P=Q^{(\mathfrak{d})}_{j,\ell,k}=A^{-j}_{(\mathfrak{d})}B_{(\mathfrak{d})}^{-\ell}(Q_0+k),
\end{equation}
where $Q_0=[0,1)^2$. The essential support of a shearlet
$\psi_{j,\ell,k}^{(\mathfrak{d})}$ in space is $P$, the volume
$\abs{P}=\abs{\text{det }A}^{-j}$ (regardless of the direction
$h,v$) and $x_P=A^{-j} B^{-\ell} k$ is the ``lower left" corner of
$P$. Let $\mathcal{Q}_{AB}:=\{Q^{(\mathfrak{d})}_{j,\ell,k}:
\mathfrak{d}=\{h,v\}, j\geq 0,\abs{\ell}\preceq 2^j,
k\in\mathbb{Z}^2\}$ and $\mathcal{Q}^{j,\ell}_{(\mathfrak{d})}:=
\{Q^{(\mathfrak{d})}_{j,\ell,k}: k\in\mathbb{Z}^d\}$. Then, for
fixed $\mathfrak{d}, j, \ell$,
$\mathcal{Q}^{j,\ell}_{(\mathfrak{d})}$ is a partition of
$\mathbb{R}^2$. Let $\mathcal{Q}^\ast_{AB}$ be the union of the indices of the shearlets and the translations of low frequency function $\varphi$. Denote $\tilde{\psi}(x)=\bar{\psi}(-x)$. We will often drop the super indices $h$ or $v$ since in makes no difference in the calculations below.
We denote for a matrix $M\in GL_2(\mathbb{R})$ the
anisotropic dilation $\varphi_M(x)=\abs{\text{det }
M}^{-1}\varphi(M^{-1}x)$. Then,
$$\psi_{A^{-j}B^{-\ell}}(x-A^{-j}B^{-\ell}k)=\abs{\text{det
}A}^{j/2}\psi_{j,\ell,k}(x)
    =\abs{P}^{-1/2}\psi_P(x),$$
and thus for $\xi\in\hat{\mathbb{R}}^2$
$$\left(\psi_{A^{-j}B^{-\ell}}(\cdot - A^{-j}B^{-\ell}k)\right)^\wedge(\xi)
    =\hat{\psi}(\xi A^{-j}B^{-\ell})\mathbf{e}^{-2\pi i\xi A^{-j}B^{-\ell}k}.$$
We also have
\begin{eqnarray}\label{e:InnrProd_equiv_Conv_Shrlts}
\nonumber
  \ip{f}{\psi_P}
    &=& \ip{f}{\psi_{j,\ell,k}} \\
\nonumber
    &=& \int_{\mathbb{R}^2} f(x) \overline{\abs{\text{det }A}^{-j/2}\psi_{A^{-j}B^{-\ell}}(x-A^{-j}B^{-\ell}k)} dx \\
    &=& \abs{P}^{1/2}(f\ast\tilde{\psi}_{A^{-j}B^{-\ell}})(x_P).
\end{eqnarray}
Let $\varphi$ be the low frequency function of the shearlet system. The
shearlet expansion for $f\in L^2(\mathbb{R}^2)$ is
\begin{eqnarray*}
  f &=& \sum_{P\in\mathcal{Q}^\ast_{AB}} \ip{f}{\psi_P}\psi_P\\
    &=& \sum_{k\in \mathbb{Z}^2} \ip{f}{\varphi_k}\varphi_k \\
    & & +
        \sum_{\mathfrak{d}=\{h,v\}}\sum_{j\geq0}\sum_{\ell=-2^j}^{2^j}\sum_{k\in\mathbb{Z}^2}
        \ip{f}{\psi^{(\mathfrak{d})}_{j,\ell,k}}\psi^{(\mathfrak{d})}_{j,\ell,k}.
\end{eqnarray*}

Define the class $\mathcal{E}$ of cartoon-like functions as those which are $C^2$ except in $C^2$ curves. It has been shown that the rate of
approximation error with $N$ shearlet coefficients $f^S_N$ is $\norm{f-f^S_N}^2_2\leq C N^{-2} (\log N)^3$, which is the theoretical optimal approximation except for the $\log$ factor. Only the curvelet system attains the same rate.


\subsection{Linogram acquisition and the pseudo-polar grid}
The PP-grid is defined by lines with equispaced slopes and intercepts with concentric squares instead of equispaced angles and concentric circles. By duality, a line through the origin in the space domain with angle $\theta$ is mapped to a line through the origin in the Fourier domain with an angle $\hat{\theta}$ orthogonal to $\theta$. Thus, the ``basically vertical" lines defining the Radon transform are mapped (via the projection slice Theorem) to ``basically horizontal" lines in the (PP) Fourier domain.

Let $m=-N/2,\ldots, N/2-1$ be the index for slopes $\frac{2m}{N}$ such that $\hat{\theta}_m=\arctan \frac{2m}{N}$ then $\hat{\theta}_m\in [-\pi/4,\pi/4)$. Let $n=-N,\ldots,N-1$ be the index for discrete distances over lines with slopes $\frac{2m}{N}$ from the center of the PP-grid. The set of points in the horizontal cone on the PP-grid is defined by $\xi_{m,n}=(\frac{n}{2N}, \frac{n}{2N}\frac{2m}{N})$. Analogously, for $m=-N/2,\ldots, N/2-1$ let $\hat{\theta}_m=\arctan \frac{2m}{N}+\pi/2$ then $\hat{\theta}_m\in [\pi/4,3\pi/4)$. The set of points in the vertical cone on the PP-grid is defined by $\xi_{m,n}=(\frac{n}{2N}\frac{2m}{N}, \frac{n}{2N})$. The pseudo-polar fast Fourier transform was developed in the fundamental paper \cite{ACDI} and was shown to be computed in $O(N^2\log N)$ flops with the Chirp-Z transform. It takes images of size $N\times N$ in the Cartesian grid to an arrange of $4N^2$ elements in the PP-grid. The oversampling factor of 4 (due to a Dirichlet interpolation kernel) allows an algebraically exact and geometrically faithful inversion both of the PP-FFT as well as the Radon transform.

The PP-grid will be used in Section \ref{S:R-inv_DscrtShrlts} when we discretize the main result.

\vskip1cm

\section{Radon inversion via shearlets}\label{S:R-inv_Shrlts_cont}
We now use results in Sections \ref{sS:RT} and \ref{S:Shrlts} to relate the properties of the Radon transform with those of the shearlet system. Let $\mathcal{R}_\theta f(s)$ be the line integral for a fixed $\theta$.

\begin{theorem}\label{t:RadonInv-Shrlts}
Let $f\in \mathcal{S}(\mathbb{R}^2)$ and $\psi$ as in Section \ref{S:Shrlts} with $\hat{\psi}_1$ even. Then, the shearlet coefficients of $f$ can be evaluated in polar coordinates at the Fourier side from its Radon transform as follows: use the Fourier slice Theorem to convert
the Radon data into polar Fourier coordinates. The coefficients can be calculated as
\begin{equation*}\label{e:ShrltCoeff-RT}
  \ip{f}{\psi_{j,\ell,k}}=\abs{P}^{1/2}\int_0^{\pi} \int_{-\infty}^\infty \hat{f}(\sigma\Theta)
        \overline{\hat{\psi}_1(2^{-2j}\sigma\cos\theta)\hat{\psi}_2(2^j\tan\theta-\ell)e^{-2\pi i \sigma\Theta x_P}} \abs{\sigma} d\sigma d\theta.
\end{equation*}
\end{theorem}
\textbf{Proof}.
From the previous notation we have that $\xi_1=\sigma\cos\theta$ and
$\xi_2=\sigma\sin\theta$. From (\ref{e:InnrProd_equiv_Conv_Shrlts}),
we apply (\ref{e:InvFrmls-RieszPotencials-special-case}) on $\tilde{\psi}_{A^{-j}B^{-\ell}}$ and
then Theorem \ref{T:BP-equiv} to obtain
\begin{eqnarray}\label{e:1st-eq-main-result}
  \nonumber
  \ip{f}{\psi_{j,\ell,k}}
    &=& \abs{P}^{1/2} (f\ast\tilde{\psi}_{A^{-j}B^{-\ell}})(x_P) \\
  \nonumber
    &=& \abs{P}^{1/2} (f\ast \frac{\mathcal{R}^\ast}{2}I^{-1} \mathcal{R}_\theta\tilde{\psi}_{A^{-j}B^{-\ell}})(x_P) \\
    &=& \abs{P}^{1/2} \frac{\mathcal{R}^\ast}{2}(\mathcal{R}_\theta
        f\ast I^{-1} \mathcal{R}_\theta\tilde{\psi}_{A^{-j}B^{-\ell}})(x_P).
\end{eqnarray}
Before proceeding we need two intermediate results. First, by
definition
\begin{eqnarray*}
  \mathcal{R}_\theta f(\ip{x_P}{\Theta}-s)
    &=& \int_{\mathbb{R}^2} f(x)\delta(\ip{x_P}{\Theta}-s-\ip{x}{\Theta}) dx \\
    &=& \int_{\mathbb{R}^2} (T_{-x_P}f)(x) \delta(-s-\ip{x}{\Theta}) dx
        = \mathcal{R}_\theta(T_{-x_P}f)(-s),
\end{eqnarray*}
and by Theorem \ref{T:Fourier-slice}
\begin{equation}\label{e:FSliceThm_shft-f}
(\mathcal{R}_\theta(T_{-x_P} f)(-\cdot))^\wedge(\sigma)=(T_{-x_P}
    f)^\wedge(-\sigma\Theta).
\end{equation}
Second, again by Theorem \ref{T:Fourier-slice}
\begin{eqnarray}\label{e:FSliceThm_shrlt}
\nonumber
    (I^{-1} \mathcal{R}_\theta\tilde{\psi}_{A^{-j}B^{-\ell}})^\wedge(\sigma)
    &=& \abs{\sigma}\mathcal{F}_{s\rightarrow\sigma} (\mathcal{R}_\theta(\tilde{\psi}_{A^{-j}B^{-\ell}}))(\sigma)\\
\nonumber
    &=& \abs{\sigma}(\tilde{\psi}_{A^{-j}B^{-\ell}})^\wedge(\sigma\Theta) \\
\nonumber
    &=& \abs{\sigma}\bar{\hat{\psi}}(\sigma\Theta A^{-j}B^{-\ell}) \\
\nonumber
    &=& \abs{\sigma}\bar{\hat{\psi}} (2^{-2j}\sigma\cos\theta, -2^{-2j}\ell\sigma\cos\theta+2^{-j}\sigma\sin\theta)\\
    &=& \abs{\sigma}\bar{\hat{\psi}}_1(2^{-2j}\sigma\cos\theta)
        \bar{\hat{\psi}}_2(2^j\tan\theta-\ell),
\end{eqnarray}
where $\hat{\psi}_1$ is an even function.

Hence, continuing with equation (\ref{e:1st-eq-main-result}) we see that Plancherel's equality, (\ref{e:FSliceThm_shft-f}) and (\ref{e:FSliceThm_shrlt}) yield
\begin{eqnarray}\label{e:chain_RT-ShrltsCoefs}
\nonumber
  \ip{f}{\psi_{j,\ell,k}}
    &=& \frac{\abs{P}^{1/2}}{2}\int_0^{2\pi} (\mathcal{R}_\theta f\ast I^{-1} \mathcal{R}_\theta\tilde{\psi}_{A^{-j}B^{-\ell}})(\ip{x_P}{\Theta}) d\theta \\
\nonumber
    &=&\frac{\abs{P}^{1/2}}{2}\int_0^{2\pi}\int_{-\infty}^\infty
        \mathcal{R}_\theta f(s)\cdot I^{-1} \mathcal{R}_\theta\tilde{\psi}_{A^{-j}B^{-\ell}}(\ip{x_P}{\Theta}-s) ds
        d\theta\\
\nonumber
    &=& \frac{\abs{P}^{1/2}}{2}\int_0^{2\pi}
        \int_{-\infty}^\infty
        \mathcal{R}_\theta f(\ip{x_P}{\Theta}-s)\cdot I^{-1} \mathcal{R}_\theta\tilde{\psi}_{A^{-j}B^{-\ell}}(s) ds d\theta \\
\nonumber
    &=& \frac{\abs{P}^{1/2}}{2}\int_0^{2\pi}
        \int_{-\infty}^\infty
        (T_{-x_P} f)^\wedge(-\sigma\Theta)\cdot \abs{\sigma} \bar{\hat{\psi}}_1(2^{-2j}\sigma\cos\theta)
        \bar{\hat{\psi}}_2(2^j\tan\theta-\ell) d\sigma d\theta\\
\nonumber
    &=& \frac{\abs{P}^{1/2}}{2}\int_0^{2\pi} \int_{-\infty}^\infty
        \hat{f}(-\sigma\Theta)e^{-2\pi i \sigma\Theta x_P} \abs{\sigma} \bar{\hat{\psi}}_1(2^{-2j}\sigma\cos\theta)
        \bar{\hat{\psi}}_2(2^j\tan\theta-\ell) d\sigma d\theta\\
\nonumber
    &=& \frac{\abs{P}^{1/2}}{2}\int_0^{2\pi} \int_{-\infty}^\infty
        \hat{f}(\sigma\Theta)e^{2\pi i \sigma\Theta x_P} \abs{\sigma} \bar{\hat{\psi}}_1(2^{-2j}\sigma\cos\theta)
        \bar{\hat{\psi}}_2(2^j\tan\theta-\ell) d\sigma d\theta\\
\nonumber
\sharp&=& \abs{P}^{1/2}\int_0^{\pi} \int_{-\infty}^\infty
        \hat{f}(\sigma\Theta)e^{2\pi i \sigma\Theta x_P} \bar{\hat{\psi}}_1(2^{-2j}\sigma\cos\theta)
        \bar{\hat{\psi}}_2(2^j\tan\theta-\ell) \abs{\sigma} d\sigma d\theta\\
\natural&=& \abs{P}^{1/2}\int_{\hat{\mathbb{R}}^2}
        \hat{f}(\xi)
        \bar{\hat{\psi}}_1(2^{-2j}\xi_1) \bar{\hat{\psi}}_2(2^j\frac{\xi_2}{\xi_1}-\ell)
        e^{2\pi i\xi x_P}d\xi,
\end{eqnarray}
where $\sharp$ is the result and, by a change of variables from polar to Cartesian,
the line marked with $\natural$ is obviously $\ip{f}{\psi_{j,\ell,k}}$ by Plancherel's equality.
\hfill $\blacksquare$ \vskip .5cm   
Theorem \ref{t:RadonInv-Shrlts} is for the shearlets associated to the horizontal cone $\mathcal{D}^h$. An analogous result holds for $\mathcal{D}^v$. For low frequencies one can choose a separable $\hat{\varphi}(\xi)=\hat{\varphi}_1(\xi_1)\hat{\varphi}_2(\xi_2)$, although we will not calculate its coefficients.

\vskip1cm

\section{The Discretization}\label{S:R-inv_DscrtShrlts}
It is line $\sharp$ in (\ref{e:chain_RT-ShrltsCoefs}) which will be discretized to match the pseudo-polar grid. Next, we briefly describe the algorithm to reconstruct a digital image of size $N\times N$ from the discrete Radon data. Assuming supp$\hat{f}\subset[-1/2,1/2]^2$, the first step is to apply a 1D Fourier transform to the $s$-variable of the Radon transform with $\theta$ fixed, i.e. from the 3rd line to the 4th line in (\ref{e:chain_RT-ShrltsCoefs}).
Every shearlet can be thought of as a ``mask" with values on the PP-grid since it is supposed to be compactly supported in the Fourier domain. The variable $\abs{\sigma}$ can be discretized as $\norm{\xi_{m,n}}$. Multiply $\hat{f}$, the shearlets, $e^{2\pi i\xi_{m,n}x_P}$ and $\norm{\xi_{m,n}}$ point-wise on the PP-grid. Finally, the shearlet coefficient is obtained by summing over slopes and distances over which the shearlet is defined on the PP-grid. This sums can be obtained by the Chirp-Z transform.

Line $\sharp$ in (\ref{e:chain_RT-ShrltsCoefs}) is discretized for the BV shearlets in time domain as follows. Then,
\begin{eqnarray}\label{e:dscrtzn_main-result}
\nonumber
  \ip{f}{\psi_{j,\ell,k}} &&  \\
    &=& \abs{P}^{1/2} \sum_{m\in\Gamma_\ell}\sum_{n\in\Delta_j} \hat{f}[m,n] \bar{\hat{\psi}}_1^{(j)}[n] \bar{\hat{\psi}}_2^{(j,\ell)}[m]
            \norm{\xi_{m,n}}e^{2\pi i\xi_{m,n}x_P},
\end{eqnarray}
where $\Gamma_\ell$ is the set of indices of lines in the PP-grid with slopes $2m/N$ (the set of $m$'s underlying the $\ell$th-directional filter or ``mask"), $\Delta_j$ is the set of indices of points from the center of the PP-grid (the set of points underlying the $j$th-level filter or ``mask"), $\hat{f}[m,n]$ is the value of $\hat{f}$ at points $\xi_{m,n}$, $\bar{\hat{\psi}}_1^{(j)}[n]$ is the value of the ``mask" $\bar{\hat{\psi}}_1^{(j)}$ at point $n$, similarly for $\bar{\hat{\psi}}_2^{(j,\ell)}[m]$, $\norm{\xi_{m,n}}=((\frac{n}{2N})^2+ (\frac{n}{2N}\frac{2m}{N})^2)^{1/2}$ and
$$\xi_{m,n}x_P=(\frac{n}{2N},\frac{n}{2N}\frac{2m}{N})(\begin{array}{c}
                                                        2^{-2j}k_1 + \ell2^{-2j}k_2 \\
                                                        2^{-j}k_2
                                                        \end{array}).$$
As mentioned earlier, the ``masks" can be stored previously. Define
$$M_{(j,\ell)}[m,n]:=\bar{\hat{\psi}}_1^{(j)}[n] \bar{\hat{\psi}}_2^{(j,\ell)}[m] \norm{\xi_{m,n}},$$
as the ``mask" at scale $j$ and direction $\ell$ with values at $\xi_{m,n}$. Define now,
$$g_{(j,\ell)}[m,n]:=\hat{f}[m,n]M_{j,\ell}[m,n].$$
Then,
\begin{eqnarray*}
  \ip{f}{\psi_{j,\ell,k}}
    && \\
    &=& \abs{P}^{1/2} \sum_{m\in\Gamma_\ell}\sum_{n\in\Delta_j} g_{(j,\ell)}[m,n]
            e^{2\pi i[\frac{n}{2N}(2^{-2j}k_1+\ell 2^{-2j}k_2) + \frac{n}{2N}\frac{2m}{N}2^{-j}k_2]}\\
    &=& \abs{P}^{1/2} \sum_{n\in\Delta_j} e^{2\pi i\frac{n}{2N}(2^{-2j}k_1+\ell 2^{-2j}k_2)}
                        \sum_{m\in\Gamma_\ell} g_{(j,\ell)}[m,n] e^{2\pi i\frac{mk_2}{N}\frac{2^{-j}n}{N}}.
\end{eqnarray*}
Denote
$$\tilde{g}_{(j,\ell)}[k_2,n]:=\sum_{m\in\Gamma_\ell} g_{(j,\ell)}[m,n] e^{2\pi i\frac{mk_2}{N}\frac{2^{-j}n}{N}}.$$
By means of the Chirp-Z transform $\tilde{g}_{(j,\ell)}[k_2,n], k_2=0,\ldots,N-1,$ can be evaluated in $O(N\log N)$ flops.
The fractional factor is $\alpha_{j,n}=\frac{2^{-j}n}{N}$. Continuing the evaluation of the shearlet coefficients we have
$$\ip{f}{\psi_{j,\ell,k}}
    =\abs{P}^{1/2} \sum_{n\in\Delta_j} e^{2\pi i\frac{n}{2N}2^{-2j}k_1} e^{2\pi i\frac{n}{2N}\ell 2^{-2j}k_2} \tilde{g}_{(j,\ell)}[k_2,n].$$
At this point, we can create a matrix $\tilde{g}_{(j,\ell)}[k_2,n], k_2=0,\ldots,N-1, n\in\Delta_j$ and for every row $k_2$ multiply entry-wise
times the vector $e^{2\pi i\frac{n}{2N}\ell 2^{-2j}k_2}$. The resulting matrix will be called $\tilde{G}_{(j,\ell)}[k_2,n]$. Finally, we evaluate
$$\ip{f}{\psi_{j,\ell,k}}
    =\abs{P}^{1/2} \sum_{n\in\Delta_j} \tilde{G}_{(j,\ell)}[k_2,n] e^{2\pi i\frac{nk_1}{N}\frac{2^{-2j}}{2}}
    =\abs{P}^{1/2} \tilde{G}_{(j,\ell)}[k_2,k_1],$$
which can be evaluated for $k_1=0,\ldots,N-1$, again, by means of the Chirp-Z transform with fractional factor $\beta_j=\frac{2^{-2j}}{2}$. Thus,
in fact, we have for fixed $j$ and $\ell$ all the shearlet coefficients $\ip{f}{\psi_{j,\ell,k}}$ for $0\leq k_1,k_2\leq N-1$.
To recover $f$ evaluate $\hat{f}=\sum \ip{f}{\psi_{j,\ell,k}} (\psi_{j,\ell,k})^\wedge$ and then apply the inverse PP-FFT.

Much of this is done by the ShearLab toolbox. More details will be given in the published paper.

Now, for fixed $\theta_m$ suppose that the vector $\mathcal{R}_{\theta_m}f(\sigma_n), n=-N,\ldots,N-1$, is additively contaminated by a noisy vector $r_m$ with i.i.d. $N(0,\sigma^2)$ entries. By the Fourier slice Theorem, the observed vector $Y_{\theta_m}=\mathcal{R}_{\theta_m}f+r_m$ turns into
$\hat{Y}_{\theta_m}=\hat{f}[m,\cdot]+\hat{r}_m$ where $\hat{r}_m$ is again a noisy vector with
i.i.d. $N(0,\sigma^2)$ entries. By linearity and Plancherel's equality we have that the error of reconstruction is
$$\norm{f-\tilde{f}}_2=\norm{\hat{f}-\hat{Y}}_2=\norm{\hat{f}-(\hat{f}+\hat{r})}_2=\norm{r}_2,$$
where $\hat{Y}$ is the non-overlapping super position of vectors $\hat{Y}_{\theta_m}$ over the PP-grid, $\tilde{f}$ is the inverse of $\hat{Y}$ and similarly for $r$.

The advantage of the proposed method over the ``more direct" PP-FFT is that, when dealing with more realistic noisy data, one can perform multiscale thresholding algorithms over the shearlet coefficients.

We now discuss the MSE over the cartoon-like functions in \cite{CaDo02} in view of our results. In \cite{CaDo02}, the inversion of the RT is obtained by the Biorthogonal Curvelet Decomposition BCD as $f = \sum_\mu [\mathcal{R}f, U_\mu] \kappa_j^{-1} \gamma_\mu$, where $\mu$ is a multi index of scale, angle and position, $\gamma_\mu$ is a curvelet, the quasi-singular value $\kappa_j^{-1}=2^j$, $U_\mu$ is one of the quasi-biorthogonal representations in the Radon domain related to the shearlet system via the intertwining relation and $[\cdot,\cdot]$ is the inner product in the Radon domain. When the RT is contaminated by Gaussian noise the empirical (or observed) coefficients are $y_\mu=[\mathcal{R}f,U_\mu]+ \epsilon[W,U_\mu]$, where $W$ is a Wiener sheet and $\epsilon$ is the noise level. From the relation $[\mathcal{R}f,U_\mu]=\kappa_j\ip{f}{\gamma_\mu}$ we can rewrite $y_\mu=\kappa_j\alpha_\mu + \epsilon \eta_\mu$, where $\alpha_\mu=\ip{f}{\gamma_\mu}$ is the curvelet coefficient of $f$ and $\eta_\mu$ is a non-i.i.d. Gaussian noise. By the tight frame property one has
$$E\norm{f-\tilde{f}}_2^2\leq E\norm{\alpha-\tilde{\alpha}}_2^2,$$
where $\tilde{f},\tilde{\alpha}$ are the estimations of the function $f$ and its its curvelets coefficients $\alpha=\{\alpha_\mu\}$, respectively. Thus, it is sufficient to estimate the coefficients $\alpha_\mu$ to estimate the MSE. This is done in \cite{CaDo02} via thresholding the observed coefficients $y_\mu$ with a threshold $t_j$ dependent on the scale $t_j=\sqrt{2\log(N_j)}\kappa_j^{-1}\epsilon$, where $N_j$ is the set of potentially non zero coefficients at scale $j$. Let $\delta(y_\mu, tj)$ be the thresholding operator. Since \cite{CaDo02} proposes
$$\tilde{\alpha}_\mu=\delta(\alpha_\mu +2^j\epsilon\eta_\mu , \sqrt{2\log (N_j)} 2^j\epsilon ),$$
this means that the estimation will only keep those $y_\mu$ greater than $\sqrt{2\log (N_j)} 2^j\epsilon$. As can be observed, the threshold is $2^j$ times the usual value $\sqrt{2\log (N_j)} \epsilon$ because the BCD multiplies noise by $2^j$. So, this estimation will discard true low-valued coefficients with high probability.

Since the analysis of our proposed inversion shows that we recover the shearlet coefficients straight from the Radon data, we obtain, at least, the same MSE $O(\epsilon^{4/5})$ as in \cite{CaDo02}. Even better, in our scheme a thresholding by $\sqrt{2\log (N_j)} \epsilon$ will recover more true shearlet coefficients, since we have $y_\mu=\alpha_\mu + \epsilon\eta_\mu$, which means that the noise level remains the same (by the orthonormality of the Fourier transform) after inversion of the RT.

\section{Simulations}\label{S:Simulations}
We show some simulations on some density functions involving rotated ellipses and squares as well as the usual Shepp-Logan head phantom.
The RT is obtained applying the PPFT to the Cartesian 2D function or image $I[i,j], i,j\in(0,\ldots,N-1)$ and mapping it to the
PP-grid by the program ppFT.m in the ShearLab toolbox. Next, as the theory mandates, we applied an 1D inverse fractional Fourier Transform and create the linogram $\mathcal{R}f[\theta_m, \cdot]$.

We add Gaussian noise to the linogram. Then, we apply the 1D fractional Fourier transform to $\mathcal{R}f[\theta_m, \cdot]$ for every angle $\theta_m$. We now have the PPFT and from here we can recover the shearlet coefficients from the ShearLab toolbox and perform the thresholding.

\begin{figure}[!t]
\centering \includegraphics[width=6in]{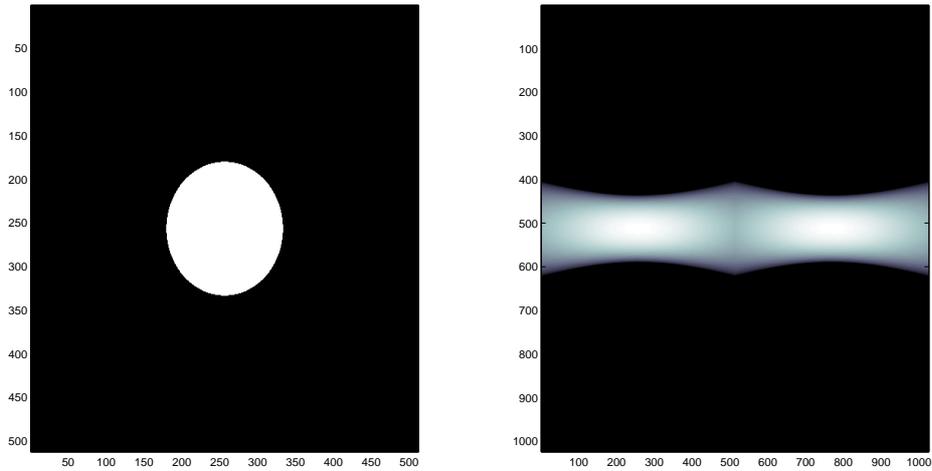} \caption{Centered circle of radius 0.3 .}\label{f:Circle3}
\end{figure}

\begin{figure}[!t]
\centering \includegraphics[width=6in]{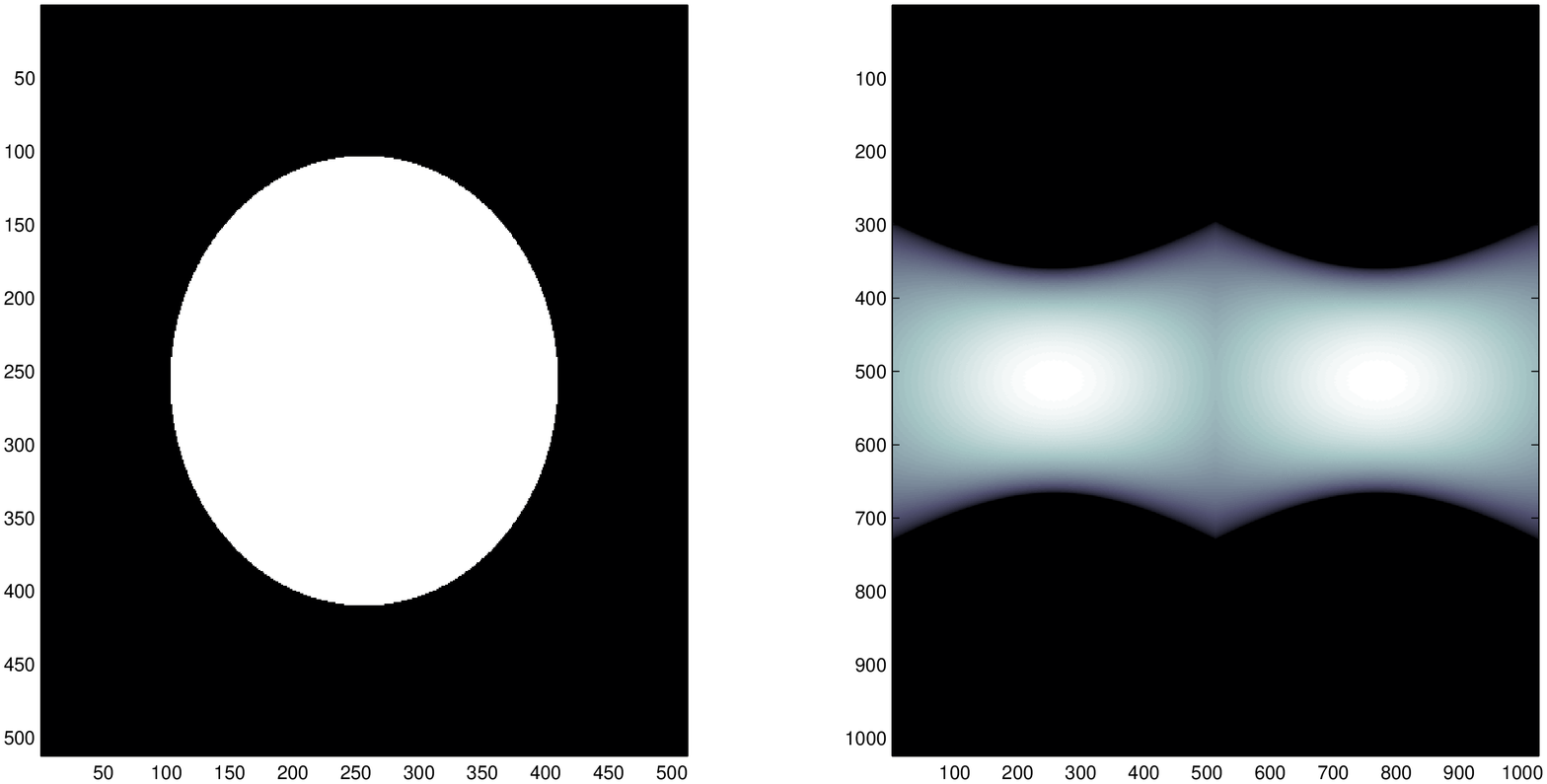} \caption{Centered circle of radius 0.6 .}\label{f:Circle6}
\end{figure}

\begin{figure}[!t]
\centering \includegraphics[width=6in]{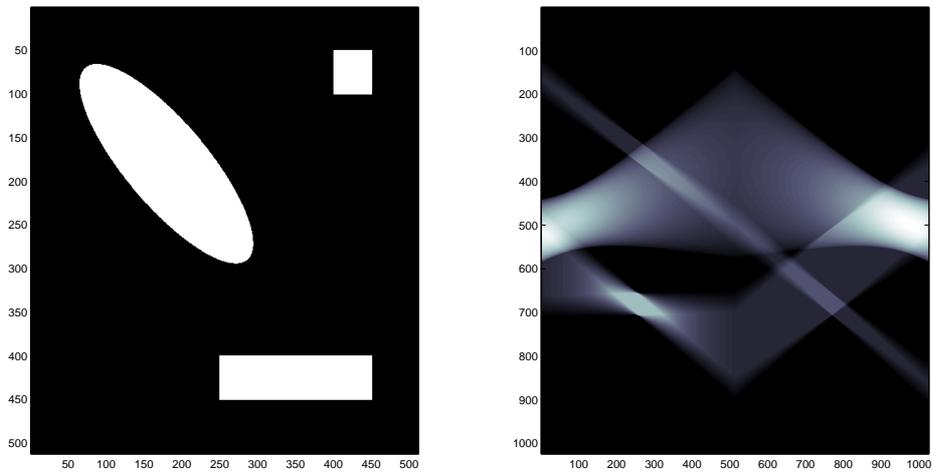} \caption{Ellipses and rectangles}\label{f:E-and-R}
\end{figure}

\begin{figure}[!t]
\centering \includegraphics[width=6in]{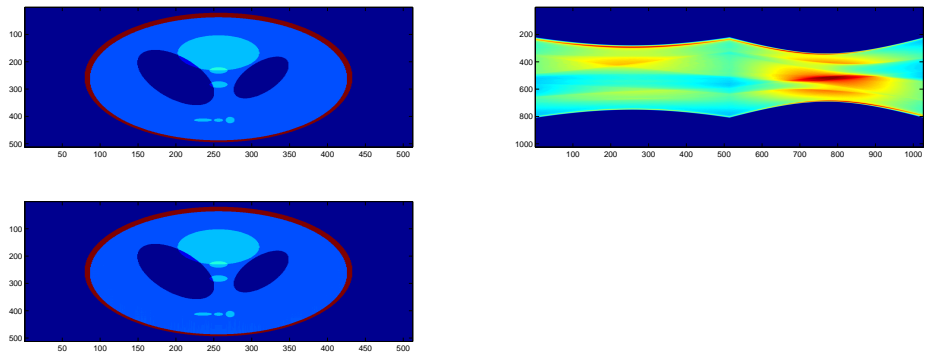} \caption{Noiseless reconstruction of Shepp-Logan phantom}\label{f:Ph-noiseless}
\end{figure}

\begin{figure}[!t]
\centering \includegraphics[width=6in]{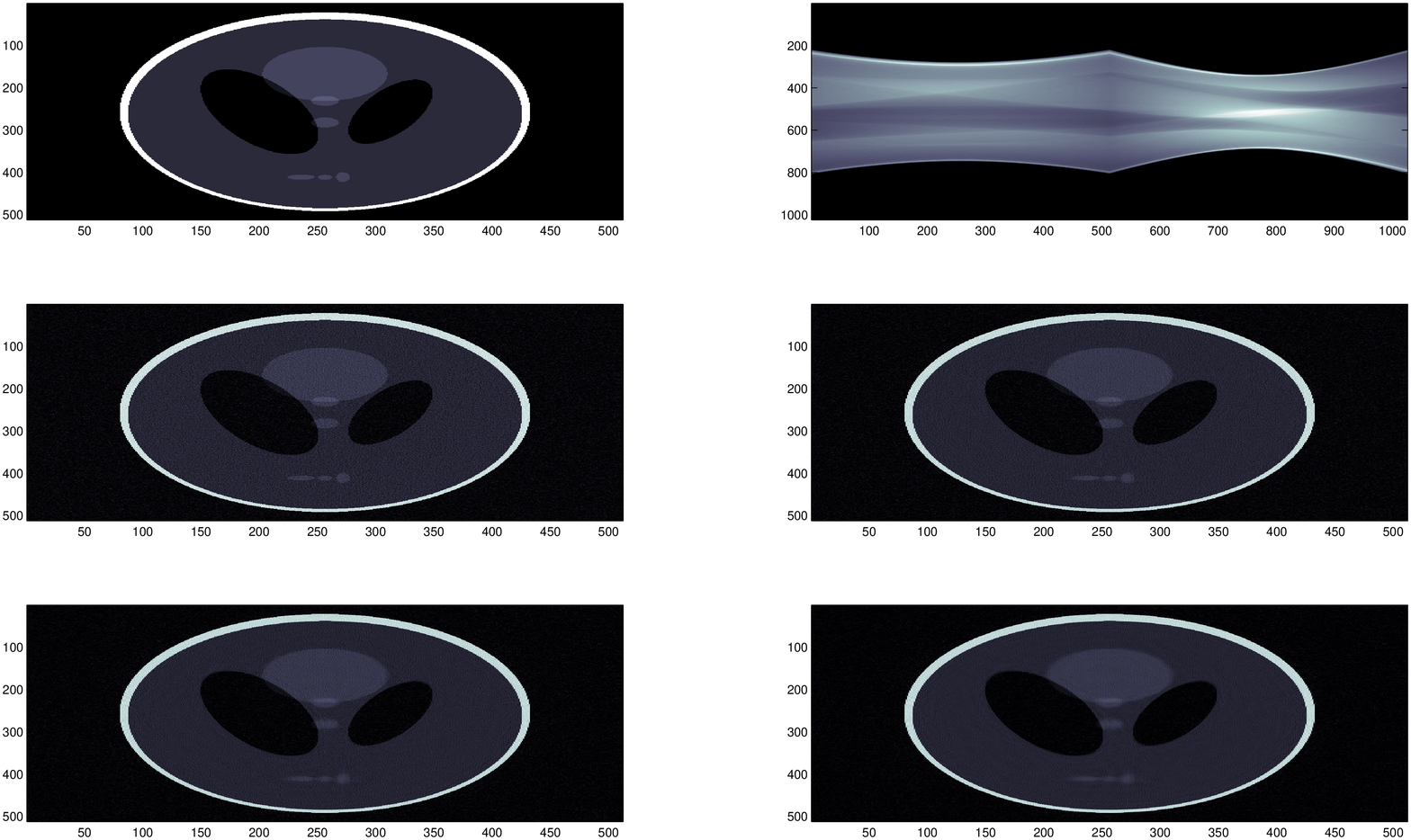} \caption{Several thresholdings}\label{f:Phs-noise}
\end{figure}



\begin{thebibliography}{99}

\parskip=0.1cm

\bibitem{ACDIS} A. Averbuch, R. Coifman, D. Donoho, M. Israeli and Y. Shkolnisky,
    ``A Framework for Discrete Integral Transformations I - The Pseudopolar Fourier Transform",
    SIAM J. Sci. Comput., 30(2), 764 - 784 (2008). \texttt{https://doi.org/10.1137/060650283}

\bibitem{BdMdVO} F. Bartolucci, F. De Mari, E. De Vito, AND F. Odone, ``Radon transform intertwines shearlets and wavelets",
http://arxiv.org/abs/1703.09578v1

\bibitem{BeWa94} C. Berenstein, D. Walnut, ``Local inversion of the Radon transform in even dimension using wavelets,"
in 75 years of Radon transform, S. Gindikin and P. Michor, Eds. Cambridge, MA: International, 1994, pp. 38-58.


\bibitem{BPTP02} S. Bonnet, F. Peyrin, F. Turjman and R. Prost, ``Multiresolution reconstruction in fan-beam tomography,"
IEEE Trans. Image Processing, vol. 11, pp. 169-176, Mar. 2002.

\bibitem{CaDo02} E. J. Candés, D.L. Donoho, Recovering edges in ill-posed
inverse problems: optimality of curvelet frames, Annals Stat., 30(3)
(2002), 784--842.

\bibitem{CEGL} F. Colonna, G. R. Easley, K. Guo, D. Labate, Radon transform
inversion using the shearlet representation, Appl. Comput. Harmon.
Anal., 29(2), p. 232-250 (2010).

\bibitem{DoVe} M. N. Do, M. Vetterli, \emph{The contourlet transform:
An efficient directional multiresolution image representation}, IEEE
Trans. Image Process. 14 (2005) 2091-2106.

\bibitem{GKL06} K. Guo, G. Kutyniok and D. Labate, \emph{Sparse multidimensional
representations using anisotropic dilation and shear operators}, in:
Wavelets and Splines, G. Chen and M. Lai (eds.), Nashboro Press,
Nashville, TN (2006), 189-201.



\bibitem{NaWu00} F. Natterer and F. Wübbeling, Mathematical Methods in
Image Reconstruction, SIAM, Philadelphia, 2000.

\bibitem{Rad17} J. Radon, Über die Bestimmung von Funktionen durch ihre Integralwerte längs gewisser Mannigfaltigkeiten,
Berichte Sächsische Akademie der Wissenschaften, Leipzig, Math-Phys. Kl., \textbf{69}, 262-267.


\end{thebibliography}
\end{document}